\documentclass[12pt]{article}
\usepackage{geometry,amsmath,amssymb,theorem,amscd,graphicx}
\usepackage {diagbox} 
\usepackage{amsmath}
\usepackage{latexsym}
\def\qed{\hfill $\Box$}

\usepackage{setspace} 
\usepackage{colortbl}
\usepackage{comment}

\newtheorem{Theorem}{Theorem}[section]

\newtheorem{Definition}[Theorem]{Definition}

\newtheorem{Corollary}[Theorem]{Corollary}

\theoremstyle{remark}
\newtheorem{Remark}{Remark}

\newcommand{\st}[2]{\genfrac{\{}{\}}{0pt}{}{#1}{#2}}

\renewcommand\mod{\mathrm{mod }}

\newcommand\pf{\noindent {\bf Proof} \quad}
\font\fivecy=wncyr5  \def\sa{{\hbox{\fivecy X}}} 
\font\sevency=wncyr7  \def\sal{{\hbox{\sevency X}}} 

\numberwithin{equation}{section}
\begin{document}
\title{Multi-indexed poly-Bernoulli numbers\large }
\author{Yuna Baba, Maki Nakasuji\footnote{The second author is supported by Grants-in-Aid for Scientific Research (C) 
22K03274.} \ and Mika Sakata
}
\date{}
\maketitle
%
%
\vskip 1cm
\par\noindent
\begin{abstract}
As properties of poly-Bernoulli numbers, a number of formulas such as the duality formula, explicit formula using the Stirling numbers of the second kind and periodicity for negative upper-index  have been established.
For the multi-indexed poly-Bernoulli numbers generalized by Kaneko-Tsumura, among such properties only the duality formula was obtained. In this paper, we restrict the double-indexed poly-Bernoulli numbers and 
show the 
explicit formula using the Stirling numbers of the second kind and periodicity for negative upper-index for them.
Further, we define the variant of multiple-indexed poly-Bernoulli numbers using the star-version of multiple-indexed logarithms and obtain the relation between this kind of double and triple-indexed poly-Bernoulli numbers with multi-indexed poly-Bernoulli numbers ahead.
\end{abstract}

{\small{Keywords: {
Bernoulli numbers}
}}
\section{Introduction}
%
%

The Bernoulli numbers $B_k$ which was defined by
$$
\sum_{i=0}^{n}\binom{n+1}{i}B_i=n+1 \quad (n=0,1,2,\ldots),
$$
were introduced by Jakob Bernoulli in his book.
One of the most famous properties of them which was conceived by Euler is the following power series expansions:
\begin{equation}\label{bernoulli}
\frac{t}{1-e^{-t}}=\sum_{n=0}^{\infty}B_n\frac{t^n}{n!}.
\end{equation}
The parameter $t$ can be real or complex.
Since their introduction, Bernoulli numbers have been extended in a various way.
poly-Bernoulli numbers by using polylogarithm series ${\rm Li}_k(z)=\sum_{m=1}^{\infty}z^m/m^k$ by Kaneko \cite{Kaneko} are the first extension:
$$\frac{{\rm Li}_k(1-e^{-t})}{1-e^{-t}}=\sum_{n=0}^{\infty}B_n^{(k)}\frac{t^n}{n!}
$$
for $k\in {\mathbb Z}$ and $n\in {\mathbb Z}_{\ge0}$ (see also Arakawa-Kaneko \cite{ArakawaKaneko1} and Arakawa-Ibukiyama-Kaneko \cite{ArakawaIbukiyamaKaneko}).
When $k=1$, $B_n^{(1)}$ is the usual Bernoulli number $B_n$.
He found the explicit formula for $B_n^{(k)}$ in terms of Stirling numbers of the second kind as a generalization of that for classical Bernoulli numbers.
\begin{Theorem}\label{explicitKaneko}(\cite[Theorem 1]{Kaneko})
For a non negative integer $n$ and an integer $k$, we have
$$B_n^{(k)}=(-1)^n\sum_{m=0}^{n}\frac{(-1)^m m! \st{n}{m}}{(m+1)^k}.
$$
Here, $\st{n}{m}$ is a Stirling number of the second kind, which is defined by the following recurrence formula:
$$\st{0}{0}=1, \st{n}{0}=\st{0}{m}=0\; (n, m\not=0), \st{n+1}{m}=\st{n}{m-1}+m\st{n}{m}.
$$
\end{Theorem}
Ohno-Sakata considered the poly-Bernoulli numbers with negative upper-index and found the following properties for their periodicity (\cite{OhnoSakata}, \cite{Sakata}).
\begin{Theorem}\label{period}(\cite[Theorem 4.2]{OhnoSakata}, \cite[Theorem 6.1]{Sakata})
Let $k, N$ be positive integers. For a prime number $p$ and positive integers $n, m\geq N$ with $n\equiv m\; (\mod\; p^{N-1}(p-1))$, we have
$$B_n^{(-k)}\equiv B_m^{(-k)} \; (\mod\; p^N).
$$
\end{Theorem}

\begin{Theorem}\label{period2}(\cite[Theorem 6.10]{Sakata})
Let $k, N$ be positive integers. For a prime number $p$ and any integer $n\geq N$, we have
$$\sum_{i=0}^{\varphi(p^N)-1}B_{n+i}^{(-k)}\equiv 0 \; (\mod\; p^N),
$$
where $\varphi(x)$ is a Euler's $\varphi$-function.
\end{Theorem}

Arakawa-Kaneko \cite{ArakawaKaneko1} introduced the multi-logarithmic function defined by
$$
{\rm Li}_{k_1, k_2, \ldots, k_r}(z)=\sum_{0<m_1<m_2<\cdots <m_r}\frac{z^{m_r}}{m_1^{k_1}m_2^{k_2}\cdots m_r^{k_r}}
$$
for $k_i\geq 1$ and $|z|<1$. And they defined a further generalization of $B_n^{(k)}$, which is known as the multiple poly-Bernoulli numbers ${\mathbb B}_n^{(k_1, \ldots, k_r)}$, as
$$
\frac{{\rm Li}_{k_1, \ldots, k_r}(1-e^{-t})}{(1-e^{-t})^r}=\sum_{n=0}^{\infty}{\mathbb B}_n^{(k_1, \ldots, k_r)}\frac{t^n}{n!}
$$
for $k_1, \ldots, k_r\in {\mathbb Z}$.
When $r=1$,  ${\mathbb B}_n^{(k_1)}=B_n^{(k_1)}$ being the poly-Bernoulli numbers mentioned as above.
Hamahata-Masubuchi and Bayad-Hamataha investigated some properties of ${\mathbb B}_n^{(k_1, \ldots, k_r)}$
(\cite{HamahataMasubuchi}, \cite{BayadHamahata}).
Kaneko-Tsumura (\cite{KanekoTsumura})
introduced further variant
using the multiple polylogarithms of $\ast$-type and of $\sal$-type defined by Goncharov (\cite{Goncharov}): 
\begin{align}
{\rm Li}_{s_1,\ldots,s_r}^\ast(z_1,\ldots,z_r)&=\sum_{1\le m_1<\cdots<m_r}\frac{z_1^{m_1}z_2^{m_2}\cdots z_r^{m_r}}{m_1^{s_1}m_2^{s_2}\cdots m_r^{s_r}}, \label{Liast}  \\
{\rm Li}_{s_1,\ldots,s_r}^\sa(z_1,\ldots,z_r)&=\sum_{1\le m_1<\cdots<m_r}\frac{z_1^{m_1}z_2^{m_2-m_1}\cdots z_r^{m_r-m_{r-1}}}{m_1^{s_1}m_2^{s_2}\cdots m_r^{s_r}} \label{List} \\
&=\sum_{\ell_1,\ldots,\ell_r=1}^{\infty}\frac{z_1^{\ell_1}z_2^{\ell_2}\cdots z_r^{\ell_r}}{\ell_1^{s_1}(\ell_1+\ell_2)^{s_2}\cdots (\ell_1+\cdots+\ell_r)^{s_r}} \notag
\end{align}
for $s_1, \ldots, s_r\in {\mathbb C}$ and $z_1, \ldots, z_r\in {\mathbb C}$.
%
They defined the multi-indexed poly-Bernoulli numbers $B_{m_1,\ldots,m_r}^{(s_1,\ldots,s_r),(d)}$ in Theorem 5.1 in \cite{KanekoTsumura}
by using multiple polylogarithms of $\sal$-type as
\begin{align}
F(x_1,\ldots,x_r;s_1,\ldots,s_r;d)&=\frac{{\rm Li}_{s_1,\ldots,s_r}^{\sa}(1-e^{-\sum_{\nu=1}^{r}x_\nu},\ldots,1-e^{-x_{r-1}-x_r},1-e^{-x_r})}{\prod_{j=1}^{d}(1-e^{-\sum_{\nu=j}^{r}x_\nu})} \notag \\
&\left(=\sum_{\ell_1,\ldots,\ell_r=1}^{\infty}\frac{\prod_{j=1}^{r}(1-e^{-\sum_{\nu=j}^{r}x_\nu})^{\ell_j-\delta_j(d)}}{\prod_{j=1}^{r}({\sum_{\nu=1}^{j}\ell_\nu})^{s_j}}\right) \notag \\
&=\sum_{m_1,\ldots,m_r=0}^{\infty}B_{m_1,\ldots,m_r}^{(s_1,\ldots,s_r),(d)}\frac{x_1^{m_1}\cdots x_r^{m_r}}{m_1!\cdots m_r!} \label{56}
\end{align}
for $x_1, \ldots, x_r\in {\mathbb C}$ and $d\in {\mathbb Z}$ ($1\leq d\leq r$).
When $d=r$, they denote $B_{m_1,\ldots,m_r}^{(s_1,\ldots,s_r),(r)}$ by ${\mathbb B}_{m_1,\ldots,m_r}^{(s_1,\ldots,s_r)}$ and proved the beautiful duality formula.
\begin{Theorem}\label{KTduality}(\cite[Theorem 5.4]{KanekoTsumura})
For $m_1, \ldots, m_r, k_1, \ldots, k_r\in {\mathbb Z}_{\geq 0}$, we have
$${\mathbb B}_{m_1,\ldots,m_r}^{(-k_1,\ldots,-k_r)}={\mathbb B}_{k_1,\ldots,k_r}^{(-m_1,\ldots,-m_r)}.
$$
\end{Theorem}
Note that for $r=1$ and $s_1=k_1\in {\mathbb Z}$, ${\mathbb B}_{m_1}^{(s_1)}={\mathbb B}_{m_1}^{(k_1)}=B_{m_1}^{(k_1)}$.
Theorem \ref{KTduality} is known as the generalization of the following duality theorem by Kaneko.
\begin{Theorem}\label{Kduality}(\cite[Theorem 2]{Kaneko})
For any $m, k\geq 0$, we have
$$B_m^{(-k)}=B_{k}^{(-m)}.
$$
\end{Theorem}
Therefore, it is natural to ask other properties such as the generalization of Theorems \ref{explicitKaneko}, \ref{period} and \ref{period2}.


On the other hand, as we can see in \eqref{Liast} and \eqref{List}, ${\rm Li}_{s_1,\ldots,s_r}^{\ast}(z_1,\ldots,z_r)$ and ${\rm Li}_{s_1,\ldots,s_r}^{\sa}(z_1,\ldots,z_r)$ are
 defined by the sum over  $m_i$ with the order without equal sign, $1\le m_1<\cdots<m_r$, so it is a natural question
 what happens in the case of equal ordering in the sum. Therefore, in this article we introduce the multiple-indexed polylogarithms. 
 \begin{Definition}
 Let
\begin{align}
{\rm Li}^{\ast, \star}_{s_1,\ldots,s_r}(z_1,\ldots,z_r)&=\sum_{1\le m_1\le \cdots\le m_r}\frac{z_1^{m_1}z_2^{m_2}\cdots z_r^{m_r}}{m_1^{s_1}m_2^{s_2}\cdots m_r^{s_r}}, \notag \\
{\rm Li}_{s_1,\ldots,s_r}^{\sa, \star}(z_1,\ldots,z_r)&=\sum_{1\le m_1\le \cdots\le m_r}\frac{z_1^{m_1}z_2^{m_2-m_1}\cdots z_r^{m_r-m_{r-1}}}{m_1^{s_1}m_2^{s_2}\cdots m_r^{s_r}} \notag \\
&=\sum_{\ell_1=1,\ell_2,\ldots,\ell_r=0}^{\infty}\frac{z_1^{\ell_1}z_2^{\ell_2}\cdots z_r^{\ell_r}}{\ell_1^{s_1}(\ell_1+\ell_2)^{s_2}\cdots (\ell_1+\cdots+\ell_r)^{s_r}}, \notag
\end{align}
for $s_1,\ldots,s_r \in \mathbb{C}$ and $z_1,\ldots,z_r\in \mathbb{C}$ with $|z_{j}|<1 \ (1\le j \le r)$.
 \end{Definition}
Similarly with \eqref{56},  we define multiple-indexed poly-Bernoulli numbers ${\mathbb B}_{m_1,\ldots,m_r}^{\star , (s_1,s_2,\ldots,s_r)}$ 
by using multiple polylogarithms of $\sal$-type as
\begin{align}
F^{\star}(x_1,\ldots,x_r;s_1,\ldots,s_r;r)
&=\frac{{\rm Li}_{s_1,\ldots,s_r}^{\sa,\star}(1-e^{-\sum_{\nu=1}^{r}x_\nu},\ldots,1-e^{-x_{r-1}-x_r},1-e^{-x_r})}{\prod_{j=1}^{r}(1-e^{-\sum_{\nu=j}^{r}x_\nu})} \notag \\
&=\sum_{m_1,\ldots,m_r=0}^{\infty}{\mathbb B}_{m_1,\ldots,m_r}^{\star,(s_1,\ldots,s_r)}\frac{x_1^{m_1}\cdots x_r^{m_r}}{m_1!\cdots m_r!}.\label{Fstar}
\end{align}
Hereafter, we may write $F^{\star}(x_1,\ldots,x_r;s_1,\ldots,s_r)$ for $F^{\star}(x_1,\ldots,x_r;s_1,\ldots,s_r; r)$. 
Then for the double-indexed and triple-indexed poly-Bernoulli numbers, we obtain the following

\begin{Theorem}\label{Thm1}
For $m_1,m_2\in\mathbb{Z}_{\ge0}$, $s_1,s_2 \in \mathbb{C}$, the double-indexed poly-Bernoulli numbers, we have
\begin{align}
m_2{\mathbb B}_{m_1,m_2-1}^{\star,(s_1,s_2)}=m_2{\mathbb B}_{m_{1},m_{2}-1}^{(s_1,s_2)}+B_{m_{2}}{\mathbb B}_{m_{1}}^{(s_1+s_2)}. \notag 
\end{align}
\end{Theorem}

\begin{Theorem}\label{Thm2}
For $m_1,m_2,m_3\in\mathbb{Z}_{\ge0}$, $s_1,s_2,s_3 \in \mathbb{C}$, the triple-indexed poly-Bernoulli numbers, we have
\begin{align}
&m_{2}m_{3}{\mathbb B}_{m_{1},m_{2}-1,m_{3}-1}^{\star,(s_1,s_2,s_3)}+m_{3}(m_{3}-1){\mathbb B}_{m_{1},m_{2},m_{3}-2}^{\star,(s_1,s_2,s_3)} \notag \\
&=m_{2}m_{3}{\mathbb B}_{m_{1},m_{2}-1,m_{3}-1}^{(s_1,s_2,s_3)}+m_{3}(m_{3}-1){\mathbb B}_{m_{1},m_{2},m_{3}-2}^{(s_1,s_2,s_3)}+m_{3}B_{m_{2}}{\mathbb B}_{m_{1},m_{3}-1}^{(s_1+s_2,s_3)} \notag \\
& \ \ \ \ \ \ \ \ +m_{2}B_{m_{3}}{\mathbb B}_{m_{1},m_{2}-1}^{(s_1,s_2+s_3)}+m_{3}B_{m_{3}-1}{\mathbb B}_{m_{1},m_{2}}^{(s_1,s_2+s_3)}+B_{m_{2}}B_{m_{3}}{\mathbb B}_{m_{1}}^{(s_1+s_2+s_3)}. \notag 
\end{align}
\end{Theorem}

In this article, in Section \ref{section2}, we discuss the generalization of Theorems \ref{explicitKaneko}, \ref{period} and \ref{period2} to the double-indexed poly-Bernoulli numbers ${\mathbb B}_{m_{1},m_{2}}^{(s_1,s_2)}$ . And in Sections \ref{section3} and \ref{section4}, we give the proof of Theorem \ref{Thm1} and \ref{Thm2}, respectively.


\section{Double-indexed poly-Bernoulli numbers}\label{section2}

In this section, we restrict double-indexed poly-Bernoulli numbers.
We first extend Theorem \ref{explicitKaneko} to double-indexed poly-Bernoulli numbers.
\begin{Theorem}\label{zenka}
For non negative integers $\ell_1,\ell_2$ and some integers $k_1,k_2$, we have
$$
\mathbb{B}_{\ell_1,\ell_2}^{(k_1,k_2)}=\sum_{0\le m_1,m_2}^{\ell_1+\ell_2}\frac{(-1)^{m_1+m_2+\ell_1+\ell_2}m_1!m_2!}{(m_1+1)^{k_1}(m_1+m_2+2)^{k_2}}\sum_{n=0}^{\ell_1+\ell_2}\st{n}{m_1}\st{\ell_1+\ell_2-n}{m_2} \binom{\ell_2}{n-\ell_1},
$$
where $\st{n}{m}$ is a Stirling number and $\binom{\ell}{n}$ is a binomial coefficient.
\end{Theorem}

\pf
This can be proved by faithful calculation. First, by the definition of ${\rm Li}^{\sa}_{k_1,k_2}$, we have 
$$
\frac{{\rm Li}^{\sa}_{k_1,k_2}(1-e^{-t_1-t_2},1-e^{-t_2})}{(1-e^{-t_1-t_2})(1-e^{-t_2})}
=\sum_{1\le m_1< m_2}\frac{(1-e^{-t_1-t_2})^{m_1-1}(1-e^{-t_2})^{m_2-m_1-1}}{m_1^{k_1}m_2^{k_2}}.$$
Putting $ m_3=m_2-m_1$, then we obtain $1\le m_3$ and 
\begin{align}
(RHS)&=\sum_{1\le m_1,m_3}\frac{(1-e^{-t_1-t_2})^{m_1-1}(1-e^{-t_2})^{m_3-1}}{m_1^{k_1}(m_1+m_3)^{k_2}} \notag \\
&=\sum_{0\le m_1,m_3}\frac{(1-e^{-t_1-t_2})^{m_1}(1-e^{-t_2})^{m_3}}{(m_1+1)^{k_1}(m_1+m_3+2)^{k_2}} \notag \\
&=\sum_{0\le m_1,m_3}\frac{(-1)^{m_1+m_3}(e^{-t_1-t_2}-1)^{m_1}(e^{-t_2}-1)^{m_3}}{(m_1+1)^{k_1}(m_1+m_3+2)^{k_2}}\frac{m_1!m_3!}{m_1!m_3!}.\notag
\end{align}
Using the properties of the Stirling numbers, we can proceed the calculation:
\begin{align}
&=\sum_{0\le m_1,m_3}\frac{(-1)^{m_1+m_3}m_1!m_3!}{(m_1+1)^{k_1}(m_1+m_3+2)^{k_2}}\sum_{n_1=0}\begin{Bmatrix} n_1 \\ m_1 \end{Bmatrix}\frac{(-t_1-t_2)^{n_1}}{n_1!}\sum_{n_2=0}\begin{Bmatrix} n_2 \\ m_3 \end{Bmatrix}\frac{(-t_2)^{n_2}}{n_2!}\notag\\
&=\sum_{0\le m_1,m_3}\frac{(-1)^{m_1+m_3}m_1!m_3!}{(m_1+1)^{k_1}(m_1+m_3+2)^{k_2}}\sum_{n_1=0}\sum_{n_2=0}^{n_1}\begin{Bmatrix} n_2 \\ m_1 \end{Bmatrix}\begin{Bmatrix} n_1-n_2 \\ m_3 \end{Bmatrix}\frac{(-t_1-t_2)^{n_2}}{n_2!}\frac{(-t_2)^{n_1-n_2}}{(n_1-n_2)!} \notag \\
&=\sum_{0\le m_1,m_3}\frac{(-1)^{m_1+m_3}m_1!m_3!}{(m_1+1)^{k_1}(m_1+m_3+2)^{k_2}}\sum_{n_1=0}\sum_{n_2=0}^{n_1}\begin{Bmatrix} n_2 \\ m_1 \end{Bmatrix}\begin{Bmatrix} n_1-n_2 \\ m_3 \end{Bmatrix}\sum_{\ell_1=0}^{n_2}\binom{n_2}{\ell_1}\frac{(-t_1)^{\ell_1}(-t_2)^{n_1-\ell_1}}{n_2!(n_1-n_2)!} \notag \\
&=\sum_{\ell_1=0}^{\infty}\sum_{0\le m_1,m_3}\frac{(-1)^{m_1+m_3}m_1!m_3!}{(m_1+1)^{k_1}(m_1+m_3+2)^{k_2}}\sum_{n_1=0}\sum_{n_2=0}^{n_1}\begin{Bmatrix} n_2 \\ m_1 \end{Bmatrix}\begin{Bmatrix} n_1-n_2 \\ m_3 \end{Bmatrix}\frac{(-t_1)^{\ell_1}(-t_2)^{n_1-\ell_1}}{\ell_1!(n_1-n_2)!(n_2-\ell_1)!} \notag \\
&=\sum_{\ell_1=0}^{\infty}\sum_{0\le m_1,m_3}\frac{(-1)^{m_1+m_3}m_1!m_3!}{(m_1+1)^{k_1}(m_1+m_3+2)^{k_2}} \notag \\
& \ \ \ \ \ \ \ \ \ \ \ \ \ \ \ \ \ \ \ \ \ \ \times\sum_{\ell_2=0}^{\infty}\sum_{n_2=0}^{\ell_1+\ell_2}\begin{Bmatrix} n_2 \\ m_1 \end{Bmatrix}\begin{Bmatrix} \ell_1+\ell_2-n_2 \\ m_3 \end{Bmatrix}\frac{(-t_1)^{\ell_1}(-t_2)^{\ell_2}}{\ell_1!(\ell_1+\ell_2-n_2)!(n_2-\ell_1)!}
\notag \\
&=\sum_{\ell_1=0}^{\infty}\sum_{0\le m_1,m_3}\frac{(-1)^{m_1+m_3}m_1!m_3!}{(m_1+1)^{k_1}(m_1+m_3+2)^{k_2}} \notag \\
& \ \ \ \ \ \ \ \ \ \ \ \ \ \ \ \ \ \ \ \ \ \ \times \sum_{\ell_2=0}^{\infty}\sum_{n_2=0}^{\ell_1+\ell_2}\begin{Bmatrix} n_2 \\ m_1 \end{Bmatrix}\begin{Bmatrix} \ell_1+\ell_2-n_2 \\ m_3 \end{Bmatrix}\binom{\ell_2}{n_2-\ell_1}\frac{(-t_1)^{\ell_1}(-t_2)^{\ell_2}}{\ell_1!\ell_2!} \notag \\
&=\sum_{\ell_1=0}^{\infty}\sum_{\ell_2=0}^{\infty}\sum_{0\le m_1,m_3}^{\ell_1+\ell_2}\frac{(-1)^{m_1+m_3+\ell_1+\ell_2}m_1!m_3!}{(m_1+1)^{k_1}(m_1+m_3+2)^{k_2}}\sum_{n_2=0}^{\ell_1+\ell_2}\begin{Bmatrix} n_2 \\ m_1 \end{Bmatrix}\begin{Bmatrix} \ell_1+\ell_2-n_2 \\ m_3 \end{Bmatrix}\binom{\ell_2}{n_2-\ell_1}\frac{t_1^{\ell_1}t_2^{\ell_2}}{\ell_1!\ell_2!}.\notag
\end{align}
Comparing the coefficients of both sides,  we obtain
$$
\mathbb{B}_{\ell_1,\ell_2}^{(k_1,k_2)}=\sum_{0\le m_1,m_3}^{\ell_1+\ell_2}\frac{(-1)^{m_1+m_3+\ell_1+\ell_2}m_1!m_3!}{(m_1+1)^{k_1}(m_1+m_3+2)^{k_2}}\sum_{n_2=0}^{\ell_1+\ell_2}\begin{Bmatrix} n_2 \\ m_1 \end{Bmatrix}\begin{Bmatrix} \ell_1+\ell_2-n_2 \\ m_3 \end{Bmatrix}\binom{\ell_2}{n_2-\ell_1}.
$$
\qed

Second, we consider the periodicity to  extend Theorem \ref{period}.
\begin{Theorem}\label{onaji}
Let $k_1,k_2,N$ be positive integers. Then 
for a prime number $p$, and positive integers $n_1,n_2,m_1,m_2\ge N$ with $n_1\equiv m_1 \ (\mod \ p^{N-1}(p-1))$
and $n_2\equiv m_2 \ (\mod \ p^{N-1}(p-1))$, we have
\begin{align}
\mathbb{B}_{n_1,n_2}^{(-k_1,-k_2)}&\equiv\mathbb{B}_{n_1,m_2}^{(-k_1,-k_2)} \ \ (\mod \ p^{N}), \label{firstas} \\
\mathbb{B}_{n_1,n_2}^{(-k_1,-k_2)}&\equiv\mathbb{B}_{m_1,n_2}^{(-k_1,-k_2)} \ \ (\mod \ p^{N}). \label{secondas}
\end{align}
\end{Theorem}

\pf
Let us begin with the first assertion.
Asuume $n_2\ge m_2$.
From Theorem $\ref{zenka}$, 
\begin{align}
&\mathbb{B}_{k_1,k_2}^{(-n_1,-n_2)}=\sum_{0\le \ell_1}^{k_1+k_2}\sum_{0\le \ell_2}^{k_1+k_2}(-1)^{\ell_1+\ell_2+k_1+k_2}\ell_1!\ell_2!(\ell_1+1)^{n_1}(\ell_1+\ell_2+2)^{n_2} \notag \\
& \ \ \ \ \ \ \ \ \ \ \ \ \ \ \ \ \ \ \ \ \ \ \ \ \ \ \times\sum_{h_2=0}^{k_1+k_2}\begin{Bmatrix} h_2 \\ \ell_1 \end{Bmatrix}\begin{Bmatrix} k_1+k_2-h_2 \\ \ell_2 \end{Bmatrix}\binom{k_2}{h_2-k_1}.
\end{align}
Dividing 
$\displaystyle{\sum_{0\le \ell_2}^{k_1+k_2}}$
into two parts 
$\displaystyle{\sum_{\begin{smallmatrix}\ell_2=0 \\ p|\ell_1+\ell_2+2 \end{smallmatrix}}^{k_1+k_2}+
\sum_{\begin{smallmatrix}\ell_2=0 \\ p\nmid\ell_1+\ell_2+2 \end{smallmatrix}}^{k_1+k_2}}
$, we obtain

\begin{align*}
\mathbb{B}_{k_1,k_2}^{(-n_1,-n_2)} 
&\equiv \sum_{0\le \ell_1}^{k_1+k_2}\sum_{\begin{smallmatrix}\ell_2=0 \\ p\nmid\ell_1+\ell_2+2 \end{smallmatrix}}^{k_1+k_2}(-1)^{\ell_1+\ell_2+k_1+k_2}\ell_1!\ell_2!(\ell_1+1)^{n_1}(\ell_1+\ell_2+2)^{n_2}\\
&\times\sum_{h_2=0}^{k_1+k_2}\begin{Bmatrix} h_2 \\ \ell_1 \end{Bmatrix}\begin{Bmatrix} k_1+k_2-h_2 \\ \ell_2 \end{Bmatrix}\binom{k_2}{h_2-k_1} \ \ \  (\mod \ p^{N}). \notag 
\end{align*}
The Euler theorem leads to
\begin{align*}
\mathbb{B}_{k_1,k_2}^{(-n_1,-n_2)}\equiv & \sum_{0\le \ell_1}^{k_1+k_2}\sum_{\begin{smallmatrix}\ell_2=0 \\ p\nmid\ell_1+\ell_2+2 \end{smallmatrix}}^{k_1+k_2}(-1)^{\ell_1+\ell_2+k_1+k_2}\ell_1!\ell_2!(\ell_1+1)^{n_1}(\ell_1+\ell_2+2)^{m_2}\\
& \times \sum_{h_2=0}^{k_1+k_2}\begin{Bmatrix} h_2 \\ \ell_1 \end{Bmatrix}\begin{Bmatrix} k_1+k_2-h_2 \\ \ell_2 \end{Bmatrix}\binom{k_2}{h_2-k_1} =\mathbb{B}_{k_1,k_2}^{(-n_1,-m_2)}\ \ \  (\mod \ p^{N}). 
\end{align*}
Applying Theorem \ref{KTduality}, we obtain
$$\mathbb{B}_{n_1,n_2}^{(-k_1,-k_2)}\equiv\mathbb{B}_{n_1,m_2}^{(-k_1,-k_2)} \ \ (\mod \ p^{N}).$$

The second assertion is obtained in a similar way. 
Assume $n_1\ge m_1$. We consider
\begin{align*}
&\mathbb{B}_{k_1,k_2}^{(-n_1,-n_2)}=\sum_{0\le \ell_2}^{k_1+k_2}\sum_{0\le \ell_1}^{k_1+k_2}(-1)^{\ell_1+\ell_2+k_1+k_2}\ell_1!\ell_2!(\ell_1+1)^{n_1}(\ell_1+\ell_2+2)^{n_2} \notag \\
& \ \ \ \ \ \ \ \ \ \ \ \ \ \ \ \ \ \ \ \ \ \ \ \ \ \ \times\sum_{h_2=0}^{k_1+k_2}\begin{Bmatrix} h_2 \\ \ell_1 \end{Bmatrix}\begin{Bmatrix} k_1+k_2-h_2 \\ \ell_2 \end{Bmatrix}\binom{k_2}{h_2-k_1}.
\end{align*}
Dividing 
$\displaystyle{\sum_{0\le \ell_1}^{k_1+k_2}}$
into two parts 
$\displaystyle{\sum_{\begin{smallmatrix}\ell_1=0 \\ p|\ell_1+1 \end{smallmatrix}}^{k_1+k_2}+
\sum_{\begin{smallmatrix}\ell_1=0 \\ p\nmid\ell_1+1 \end{smallmatrix}}^{k_1+k_2}}
$, we obtain

\begin{align*}
\mathbb{B}_{k_1,k_2}^{(-n_1,-n_2)}&\equiv \sum_{0\le \ell_2}^{k_1+k_2}\sum_{\begin{smallmatrix}\ell_1=0 \\ p\nmid\ell_1+1 \end{smallmatrix}}^{k_1+k_2}(-1)^{\ell_1+\ell_2+k_1+k_2}\ell_1!\ell_2!(\ell_1+1)^{n_1}(\ell_1+\ell_2+2)^{n_2}\\
&\times \sum_{h_2=0}^{k_1+k_2}\begin{Bmatrix} h_2 \\ \ell_1 \end{Bmatrix}\begin{Bmatrix} k_1+k_2-h_2 \\ \ell_2 \end{Bmatrix}\binom{k_2}{h_2-k_1} \ \ \  (\mod \ p^{N}). \notag 
\end{align*}
The Euler theorem leads to
\begin{align*}
\mathbb{B}_{k_1,k_2}^{(-n_1,-n_2)}\equiv \sum_{0\le \ell_2}^{k_1+k_2}\sum_{\begin{smallmatrix}\ell_1=0 \\ p\nmid\ell_1+1 \end{smallmatrix}}^{k_1+k_2}(-1)^{\ell_1+\ell_2+k_1+k_2}\ell_1!\ell_2!(\ell_1+1)^{m_1}(\ell_1+\ell_2+2)^{n_2}\\
\times \sum_{h_2=0}^{k_1+k_2}\begin{Bmatrix} h_2 \\ \ell_1 \end{Bmatrix}\begin{Bmatrix} k_1+k_2-h_2 \\ \ell_2 \end{Bmatrix}\binom{k_2}{h_2-k_1} =\mathbb{B}_{k_1,k_2}^{(-m_1,-n_2)}\ \ \  (\mod \ p^{N}). \notag
\end{align*}
Applying Theorem \ref{KTduality}, we obtain the second assertion:
$$\mathbb{B}_{n_1,n_2}^{(-k_1,-k_2)}\equiv\mathbb{B}_{m_1,n_2}^{(-k_1,-k_2)} \ \ (\mod \ p^{N}).$$
\qed

Next, we consider the periodicity to  extend Theorem \ref{period2}.

\begin{Theorem}\label{0ninaru}
Let $k_1,k_2,N$ be positive integers. 
For a prime number $p$, and integers $n_1,n_2\ge N$, we have
\begin{align}
\sum_{i=0}^{\varphi(p^{N})-1}\mathbb{B}_{n_1+i,n_2}^{(-k_1,-k_2)}&\equiv0 \ \ (\mod \ p^{N}), \label{B=0n1} \\
\sum_{i=0}^{\varphi(p^{N})-1}\mathbb{B}_{n_1,n_2+i}^{(-k_1,-k_2)}&\equiv0 \ \ (\mod \ p^{N}). \label{B=0n2}
\end{align}
\end{Theorem}

\pf
For \eqref{B=0n1}, we first consider the case of $n_1=N$.
From Theorems \ref{KTduality} and \ref{zenka},
\begin{align*}
&\sum_{i=0}^{\varphi(p^{N})-1}\mathbb{B}_{N+i,n_2}^{(-k_1,-k_2)} 
=\sum_{i=0}^{\varphi(p^{N})-1}\mathbb{B}_{k_1,k_2}^{(-N-i,-n_2)} \notag \\
&=\sum_{i=0}^{\varphi(p^{N})-1}\sum_{\ell_1,\ell_2=0}^{k_1+k_2}(-1)^{\ell_1+\ell_2+k_1+k_2}\ell_1!\ell_2!(\ell_1+1)^{N+i}(\ell_1+\ell_2+2)^{n_2}\\
& \ \ \ \ \ \ \ \ \ \ \ \ \ \ \ \ \ \ \ \ \ \times \sum_{h_2=0}^{k_1+k_2}\begin{Bmatrix} h_2 \\ \ell_1 \end{Bmatrix}\begin{Bmatrix} k_1+k_2-h_2 \\ \ell_2 \end{Bmatrix}\binom{k_2}{h_2-k_1} \notag \\
&=\sum_{\ell_1,\ell_2=0}^{k_1+k_2}(-1)^{\ell_1+\ell_2+k_1+k_2}\ell_1!\ell_2!(\ell_1+\ell_2+2)^{n_2}\\
& \ \ \ \ \ \ \ \ \ \ \ \ \ \ \ \ \ \ \ \ \ \times \sum_{h_2=0}^{k_1+k_2}\begin{Bmatrix} h_2 \\ \ell_1 \end{Bmatrix}\begin{Bmatrix} k_1+k_2-h_2 \\ \ell_2 \end{Bmatrix}\binom{k_2}{h_2-k_1}\frac{(\ell_1+1)^{N}\{(\ell_1+1)^{\varphi(p^{N})}-1\}}{\ell_1}.\notag 
\end{align*}
The right-hand side above identically $0$ modulo $p^{N}$ from the Euler theorem and we get the assertion.
In case of $n_1>N$, applying the case of $n_1=N$ and Theorem \ref{onaji} leads to the desired result.
For \eqref{B=0n2}, we similarly consider the case of $n_2=N$ first.
From Theorems \ref{KTduality} and \ref{zenka},
\begin{align*}
&\sum_{i=0}^{\varphi(p^{N})-1}\mathbb{B}_{n,_1N+i}^{(-k_1,-k_2)} 
=\sum_{i=0}^{\varphi(p^{N})-1}\mathbb{B}_{k_1,k_2}^{(-n_1,-N-i)} \notag \\
&=\sum_{i=0}^{\varphi(p^{N})-1}\sum_{\ell_1,\ell_2=0}^{k_1+k_2}(-1)^{\ell_1+\ell_2+k_1+k_2}\ell_1!\ell_2!(\ell_1+1)^{n_1}(\ell_1+\ell_2+2)^{N+i}\\
& \ \ \ \ \ \ \ \ \ \ \ \ \ \ \ \ \ \ \ \ \ \times\sum_{h_2=0}^{k_1+k_2}\begin{Bmatrix} h_2 \\ \ell_1 \end{Bmatrix}\begin{Bmatrix} k_1+k_2-h_2 \\ \ell_2 \end{Bmatrix}\binom{k_2}{h_2-k_1} \notag \\
&=\sum_{\ell_1,\ell_2=0}^{k_1+k_2}(-1)^{\ell_1+\ell_2+k_1+k_2}\ell_1!\ell_2!(\ell_1+1)^{n_1}\\
& \ \ \ \ \ \ \ \ \ \ \ \ \times \sum_{h_2=0}^{k_1+k_2}\begin{Bmatrix} h_2 \\ \ell_1 \end{Bmatrix}\begin{Bmatrix} k_1+k_2-h_2 \\ \ell_2 \end{Bmatrix}\binom{k_2}{h_2-k_1}\frac{(\ell_1+\ell_2+2)^{N}\{(\ell_1+\ell_2+2)^{\varphi(p^{N})}-1\}}{\ell_1+\ell_2+1}.
\end{align*}
The right-hand side above also identically $0$ modulo $p^{N}$ from the Euler theorem and we get the assertion.
The case of $n_2>N$ is similarly derived.
\qed

\begin{Corollary}
Let $k_1,k_2,M$ be positive integers.
For $M=p_1^{e_1}p_2^{e_2}\cdots p_m^{e_m}$, where $p_i$'s are relatively prime, and any integer $n\ge e_1,e_2,\ldots,e_m$, we have
\begin{align}
\sum_{i=0}^{\varphi(M)-1}\mathbb{B}_{n_1+i,n_2}^{(-k_1,-k_2)}&\equiv0 \ \ (\mod \ M), \label{Cor1} \\
\sum_{i=0}^{\varphi(M)-1}\mathbb{B}_{n_1,n_2+i}^{(-k_1,-k_2)}&\equiv0 \ \ (\mod \ M). \label{Cor2}
\end{align}
\end{Corollary}

\pf
From Theorem \ref{onaji}, for $1\le r\le m$, we have
\begin{align}
\sum_{i=0}^{\varphi(M)-1}\mathbb{B}_{n_1+i,n_2}^{(-k_1,-k_2)}&=\sum_{i=0}^{\varphi(p_1^{e_1})\cdots\varphi(p_m^{e_m})-1}\mathbb{B}_{n_1+i,n_2}^{(-k_1,-k_2)} \notag \\
&\equiv\prod_{\begin{smallmatrix}j_1=0 \\ j\neq r \end{smallmatrix}}^{m}\varphi(p_j^{e_j})\sum_{i=0}^{\varphi(p_r^{e_r})-1}\mathbb{B}_{n_1+i,n_2}^{(-k_1,-k_2)} \ \ (\mod \ p_r^{e_r}). \notag
\end{align}
The right-hand side above also identically $0$ modulo $p_r^{e_r}$ from Theorem \ref{0ninaru}, we can get \eqref{Cor1}.
We similarly consider the case of $n_2$ and obtain \eqref{Cor2}.
\qed

\begin{Remark}
As we can see, Theorem \ref{zenka} is the key theorem for Theorems \ref{onaji} and \ref{0ninaru}.
Therefore, if 
Theorem \ref{zenka} can be generalized for $n$-indexed poly-Bernoulli numbers with $n\geq 3$, 
the related periodicity theorems 
may also be derived in a similar way.
However, we have not pursued here because the binomial coefficients in these cases are complicated.  
\end{Remark}

\begin{Remark}
According to our calculations (by using Mathematica), 
 there seems to have finer periodicities in \eqref{firstas} and \eqref{secondas} sometimes, but not always.
This means that we may be able to obtain better results under suitable conditions.
\end{Remark}

\section{Double-indexed star-poly-Bernoulli numbers}\label{section3}
This section is devoted to the proof of Theorem \ref{Thm1}, where
we restrict double-indexed polylogarithms.
By the harmonic calculations, we have 
%
\begin{align}
{\rm Li}_{k_1,k_2}^{\sa,\star}(z_1,z_2) 
&=\sum_{1\le m_1\le m_2}\frac{z_1^{m_1}z_2^{m_2-m_1}}{m_1^{k_1}m_2^{k_2}} \notag \\
&=\sum_{1\le m_1<m_2}\frac{z_1^{m_1}z_2^{m_2-m_1}}{m_1^{k_1}m_2^{k_2}}+\sum_{1\le m_1}\frac{z_1^{m_1}}{m_1^{k_1+k_2}} \notag \\
&={\rm Li}_{k_1,k_2}^{\sa}(z_1,z_2)+{\rm Li}_{k_1+k_2}^{\sa}(z_1). \notag 
\end{align}
So, $F^{\star}(x_1,x_2;s_1,s_2)$ defined in \eqref{Fstar} is written as
\begin{align}
F^{\star}(x_1,x_2;s_1,s_2) 
&=\frac{{\rm Li}_{s_1,s_2}^{\sa,\star}(1-e^{-x_1-x_2},1-e^{-x_2})}{(1-e^{-x_1-x_2})(1-e^{-x_2})} \notag \\
&=\frac{{\rm Li}_{s_1,s_2}^{\sa}(1-e^{-x_1-x_2},1-e^{-x_2})}{(1-e^{-x_1-x_2})(1-e^{-x_2})}+\frac{1}{(1-e^{-x_2})}\frac{{\rm Li}_{s_1+s_2}^{\sa}(1-e^{-x_1-x_2})}{(1-e^{-x_1-x_2})}. \notag \\
\end{align}
Using \eqref{56}, this can be expressed in terms of double-indexed poly-Bernoulli numbers,
\begin{align}
F^{\star}(x_1,x_2;s_1,s_2) &=\sum_{m_{11},m_{12}=0}^{\infty}{\mathbb B}_{m_{11},m_{12}}^{(s_1,s_2)}\frac{x_1^{m_{11}}x_2^{m_{12}}}{m_{11}!m_{12}!}+\frac{1}{(1-e^{-x_2})}\sum_{m_{21}=0}^{\infty}{\mathbb B}_{m_{21}}^{(s_1+s_2)}\frac{(x_1+x_2)^{m_{21}}}{m_{21}!}.
\end{align}
Further, by \eqref{bernoulli}, we have
\begin{align}
&F^{\star}(x_1,x_2;s_1,s_2) \notag \\
&=\sum_{m_{11},m_{12}=0}^{\infty}{\mathbb B}_{m_{11},m_{12}}^{(s_1,s_2)}\frac{x_1^{m_{11}}x_2^{m_{12}}}{m_{11}!m_{12}!}+\frac{1}{x_2}\sum_{{n_{2}}=0}^{\infty}B_{n_{2}}\frac{x_2^{n_{2}}}{{n_{2}}!}\sum_{m_{21}=0}^{\infty}{\mathbb B}_{m_{21}}^{(s_1+s_2)}\frac{(x_1+x_2)^{m_{21}}}{m_{21}!} \notag \\
&=\sum_{m_{11},m_{12}=0}^{\infty}{\mathbb B}_{m_{11},m_{12}}^{(s_1,s_2)}\frac{x_1^{m_{11}}x_2^{m_{12}}}{m_{11}!m_{12}!}+\frac{1}{x_2}\sum_{{n_{2}}=0}^{\infty}B_{n_{2}}\frac{x_2^{n_{2}}}{{n_{2}}!}\sum_{m_{21}=0}^{\infty}{\mathbb B}_{m_{21}}^{(s_1+s_2)}\frac{\sum_{i_{2}=0}^{m_{21}}\binom{m_{21}}{i_{2}}x_1^{i_{2}}x_2^{m_{21}-i_{2}}}{m_{21}!} \notag \\
&=\sum_{m_{11},m_{12}=0}^{\infty}{\mathbb B}_{m_{11},m_{12}}^{(s_1,s_2)}\frac{x_1^{m_{11}}x_2^{m_{12}}}{m_{11}!m_{12}!}+\frac{1}{x_2}\sum_{n_{2},m_{21}=0}^{\infty}B_{n_{2}}{\mathbb B}_{m_{21}}^{(s_1+s_2)}\sum_{i_{2}=0}^{m_{21}}\binom{m_{21}}{i_{2}}\frac{x_1^{i_{2}}x_2^{m_{21}-i_{2}+n_{2}}}{n_{2}!m_{21}!} \notag \\
&=\sum_{m_{1},m_{2}=0}^{\infty}{\mathbb B}_{m_{1},m_{2}}^{(s_1,s_2)}\frac{x_1^{m_{1}}x_2^{m_{2}}}{m_{1}!m_{2}!}+\frac{1}{x_2}\sum_{m_{1},m_{2}=0}^{\infty}B_{m_{2}}{\mathbb B}_{m_{1}}^{(s_1+s_2)}\sum_{i_{2}=0}^{m_{1}}\binom{m_{1}}{i_{2}}\frac{x_1^{i_{2}}x_2^{m_{1}-i_{2}+m_{2}}}{m_{1}!m_{2}!} \notag 
\end{align}
Therefore, from the definition of double-indexed poly-Bernoulli numbers, we have
\begin{align}
&x_2\sum_{m_1,m_2=0}^{\infty}{\mathbb B}_{m_1,m_2}^{\star,(s_1,s_2)}\frac{x_1^{m_1}x_2^{m_2}}{m_1!m_2!} =x_2F^{\star}(x_1,x_2;s_1,s_2) \notag \\
&=x_2\sum_{m_{1},m_{2}=0}^{\infty}{\mathbb B}_{m_{1},m_{2}}^{(s_1,s_2)}\frac{x_1^{m_{1}}x_2^{m_{2}}}{m_{1}!m_{2}!}+\sum_{m_{1},m_{2}=0}^{\infty}B_{m_{2}}{\mathbb B}_{m_{1}}^{(s_1+s_2)}\sum_{i_{2}=0}^{m_{1}}\binom{m_{1}}{i_{2}}\frac{x_1^{i_{2}}x_2^{m_{1}-i_{2}+m_{2}}}{m_{1}!m_{2}!} \notag \\
&=\sum_{m_{1},m_{2}=0}^{\infty}m_2{\mathbb B}_{m_{1},m_{2}-1}^{(s_1,s_2)}\frac{x_1^{m_{1}}x_2^{m_{2}}}{m_{1}!m_{2}!}+\sum_{m_{1},m_{2}=0}^{\infty}B_{m_{2}}{\mathbb B}_{m_{1}}^{(s_1+s_2)}\sum_{i_{2}=0}^{m_{1}}\binom{m_{1}}{i_{2}}\frac{x_1^{i_{2}}x_2^{m_{1}-i_{2}+m_{2}}}{m_{1}!m_{2}!}.\notag \\
\end{align}
Theorem\ref{Thm1} is obtained by comparing the coefficients of this equation.

Since for an odd number $n_2$ where $n_2>1$ $B_{n_{2}}=0$, Theorem \ref{Thm1} becomes
$$
{\mathbb B}_{n_1,n_2-1}^{\star,(s_1,s_2)}={\mathbb B}_{n_{1},n_{2}-1}^{(s_1,s_2)}.
$$
This gives the following corollary.

\begin{Corollary}
Let $k_1,k_2,N$ be positive integers and $n_2$ be an odd number where $n_2>1$. 
Then for a prime number $p$, and integers $n_1,n_2,m_1,m_2\ge N$ with $n_1\equiv m_1 \ (\mod \ p^{N-1}(p-1))$
and $n_2\equiv m_2 \ (\mod \ p^{N-1}(p-1))$, we have
\begin{align}
\mathbb{B}_{n_1,n_2-1}^{\star,(-k_1,-k_2)}&\equiv\mathbb{B}_{n_1,m_2-1}^{(-k_1,-k_2)} \ \ (\mod \ p^{N}), \quad
\mathbb{B}_{n_1,n_2-1}^{\star,(-k_1,-k_2)}\equiv\mathbb{B}_{m_1,n_2-1}^{(-k_1,-k_2)} \ \ (\mod \ p^{N}), \notag
\end{align}
\begin{align}
\sum_{i=0}^{\varphi(p^{N})-1}\mathbb{B}_{n_1+i,n_2-1}^{\star,(-k_1,-k_2)}&\equiv0 \ \ (\mod \ p^{N}), \quad
\sum_{i=0}^{\varphi(p^{N})-1}\mathbb{B}_{n_1,n_2-1+i}^{\star,(-k_1,-k_2)}\equiv0 \ \ (\mod \ p^{N}). \notag
\end{align}
Let $M$ be positive integers.
For $M=p_1^{e_1}p_2^{e_2}\cdots p_m^{e_m}$, where $p_i$'s are relatively prime, and any integer $n\ge e_1,e_2,\ldots,e_m$, we have
\begin{align}
\sum_{i=0}^{\varphi(M)-1}\mathbb{B}_{n_1+i,n_2-1}^{\star,(-k_1,-k_2)}&\equiv0 \ \ (\mod \ M), \quad
\sum_{i=0}^{\varphi(M)-1}\mathbb{B}_{n_1,n_2-1+i}^{\star,(-k_1,-k_2)}\equiv0 \ \ (\mod \ M). \notag
\end{align}
\end{Corollary}


\section{Triple-indexed poly-Bernoulli numbers}\label{section4}
In this section, we restrict triple-indexed polylogarithms and give similar discussion in the previous section to lead the proof of Theorem \ref{Thm2}.
By the harmonic calculations, we have 
%
%
\begin{align}
&{\rm Li}_{k_1,k_2,k_3}^{\sa,\star}(z_1,z_2,z_3) \notag \\
&=\sum_{1\le m_1\le m_2\le m_3}\frac{z_1^{m_1}z_2^{m_2-m_1}z_3^{m_3-m_2}}{m_1^{k_1}m_2^{k_2}m_3^{k_3}} \notag \\
&=\sum_{1\le m_1<m_2<m_3}\frac{z_1^{m_1}z_2^{m_2-m_1}z_3^{m_3-m_2}}{m_1^{k_1}m_2^{k_2}m_3^{k_3}}+\sum_{1\le m_1<m_3}\frac{z_1^{m_1}z_3^{m_3-m_1}}{m_1^{k_1+k_2}m_3^{k_3}}\notag\\
&\hspace{3mm}+\sum_{1\le m_1<m_2}\frac{z_1^{m_1}z_2^{m_2-m_1}}{m_1^{k_1}m_2^{k_2+k_3}}+\sum_{1\le m_1}\frac{z_1^{m_1}}{m_1^{k_1+k_2+k_3}} \notag \\
&={\rm Li}_{k_1,k_2,k_3}^{\sa}(z_1,z_2,z_3)+{\rm Li}_{k_1+k_2,k_3}^{\sa}(z_1,z_3)+{\rm Li}_{k_1,k_2+k_3}^{\sa}(z_1,z_2)+{\rm Li}_{k_1+k_2+k_3}^{\sa}(z_1) .\notag 
\end{align}
$F^{\star}(x_1,x_2,x_3;s_1,s_2,s_3)$ is written as

\begin{align*}
F^{\star}(x_1,x_2,x_3;s_1,s_2,s_3) 
&=\frac{{\rm Li}_{s_1,s_2,s_3}^{\sa,\star}(1-e^{-x_1-x_2-x_3},1-e^{-x_2-x_3},1-e^{-x_3})}{(1-e^{-x_1-x_2-x_3})(1-e^{-x_2-x_3})(1-e^{-x_3})} 
\\
&=\frac{{\rm Li}_{s_1,s_2,s_3}^{\sa}(1-e^{-x_1-x_2-x_3},1-e^{-x_2-x_3},1-e^{-x_3})}{(1-e^{-x_1-x_2-x_3})(1-e^{-x_2-x_3})(1-e^{-x_3})}\\
& \quad \quad +\frac{1}{(1-e^{-x_2-x_3})}\frac{{\rm Li}_{s_1+s_2,s_3}^{\sa}(1-e^{-x_1-x_2-x_3},1-e^{-x_3})}{(1-e^{-x_1-x_2-x_3})(1-e^{-x_3})} \\
&  \quad \quad  +\frac{1}{(1-e^{-x_3})}\frac{{\rm Li}_{s_1,s_2+s_3}^{\sa}(1-e^{-x_1-x_2-x_3},1-e^{-x_2-x_3})}{(1-e^{-x_1-x_2-x_3})(1-e^{-x_2-x_3})}\\
& \quad \quad +\frac{1}{(1-e^{-x_2-x_3})(1-e^{-x_3})}\frac{{\rm Li}_{s_1+s_2+s_3}^{\sa}(1-e^{-x_1-x_2-x_3})}{(1-e^{-x_1-x_2-x_3})}. \notag
\end{align*}
In terms of triple-indexed poly-Bernoulli numbers and classical Bernoulli numbers, it is

\begin{align*}
&F^{\star}(x_1,x_2,x_3;s_1,s_2,s_3) \notag \\
&=\sum_{m_{11},m_{12},m_{13}=0}^{\infty}{\mathbb B}_{m_{11},m_{12},m_{13}}^{(s_1,s_2,s_3)}\frac{x_1^{m_{11}}x_2^{m_{12}}x_3^{m_{13}}}{m_{11}!m_{12}!m_{13}!}+\frac{1}{(1-e^{-x_2-x_3})}\sum_{m_{21},m_{22}=0}^{\infty}{\mathbb B}_{m_{21},m_{22}}^{(s_1+s_2,s_3)}\frac{(x_1+x_2)^{m_{21}}x_3^{m_{22}}}{m_{21}!m_{22}!} \notag \\
& \quad \quad  +\frac{1}{(1-e^{-x_3})}\sum_{m_{31},m_{32}=0}^{\infty}{\mathbb B}_{m_{31},m_{32}}^{(s_1,s_2+s_3)}\frac{x_1^{m_{31}}(x_2+x_3)^{m_{32}}}{m_{31}!m_{32}!}\\
&  \quad \quad +\frac{1}{(1-e^{-x_2-x_3})(1-e^{-x_3})}\sum_{m_{41}=0}^{\infty}{\mathbb B}_{m_{41}}^{(s_1+s_2+s_3)}\frac{(x_1+x_2+x_3)^{m_{41}}}{m_{41}!} \notag \\
&=\sum_{m_{11},m_{12},m_{13}=0}^{\infty}{\mathbb B}_{m_{11},m_{12},m_{13}}^{(s_1,s_2,s_3)}\frac{x_1^{m_{11}}x_2^{m_{12}}x_3^{m_{13}}}{m_{11}!m_{12}!m_{13}!}\\
&  \quad \quad +\frac{1}{(x_2+x_3)}\sum_{{n_{2}}=0}^{\infty}{B}_{n_{2}}\frac{(x_2+x_3)^{n_{2}}}{{n_{2}}!}\sum_{m_{21},m_{22}=0}^{\infty}{\mathbb B}_{m_{21},m_{22}}^{(s_1+s_2,s_3)}\frac{(x_1+x_2)^{m_{21}}x_3^{m_{22}}}{m_{21}!m_{22}!} \notag \\
& \quad \quad  +\frac{1}{x_3}\sum_{{n_{3}}=0}^{\infty}B_{n_{3}}\frac{x_3^{n_{3}}}{{n_{3}}!}\sum_{m_{31},m_{32}=0}^{\infty}{\mathbb B}_{m_{31},m_{32}}^{(s_1,s_2+s_3)}\frac{x_1^{m_{31}}(x_2+x_3)^{m_{32}}}{m_{31}!m_{32}!} \notag \\
& \quad \quad  +\frac{1}{(x_2+x_3)x_3}\sum_{{n_{41}}=0}^{\infty}B_{n_{41}}\frac{(x_2+x_3)^{n_{41}}}{{n_{41}}!}\sum_{{n_{42}}=0}^{\infty}B_{n_{42}}\frac{x_3^{n_{42}}}{{n_{42}}!}\sum_{m_{41}=0}^{\infty}{\mathbb B}_{m_{41}}^{(s_1+s_2+s_3)}\frac{(x_1+x_2+x_3)^{m_{41}}}{m_{41}!} 
\end{align*}
\begin{align*}
&=\sum_{m_{1},m_{2},m_{3}=0}^{\infty}{\mathbb B}_{m_{1},m_{2},m_{3}}^{(s_1,s_2,s_3)}\frac{x_1^{m_{1}}x_2^{m_{2}}x_3^{m_{3}}}{m_{1}!m_{2}!m_{3}!} \notag \\
& \ \ \ \ \ \ \ \ +\frac{1}{(x_2+x_3)}\sum_{m_{2},m_{1},m_{3}=0}^{\infty}B_{m_{2}}{\mathbb B}_{m_{1},m_{3}}^{(s_1+s_2,s_3)}\sum_{i_{2}=0}^{m_{2}}\binom{m_{2}}{i_{2}}\sum_{j_{21}=0}^{m_{1}}\binom{m_{1}}{j_{21}}\frac{x_1^{j_{21}}x_2^{i_{2}+m_{1}-j_{21}}x_3^{m_{2}-i_{2}+m_{3}}}{m_{1}!m_{2}!m_{3}!} \notag \\
& \ \ \ \ \ \ \ \ +\frac{1}{x_3}\sum_{m_{3},m_{1},m_{2}=0}^{\infty}B_{m_{3}}{\mathbb B}_{m_{1},m_{2}}^{(s_1,s_2+s_3)}\sum_{j_{32}=0}^{m_{2}}\binom{m_{2}}{j_{32}}\frac{x_1^{m_{1}}x_2^{j_{32}}x_3^{m_{3}+m_{2}-j_{32}}}{m_{1}!m_{2}!m_{3}!} \notag \\
& \ \ \ \ \ \ \ \ +\frac{1}{(x_2+x_3)x_3}\sum_{m_{2},m_{3},m_{1}=0}^{\infty}B_{m_{2}}{B}_{m_{3}}{\mathbb B}_{m_{1}}^{(s_1+s_2+s_3)} \notag \\
& \ \ \ \ \ \ \ \ \ \ \ \ \ \ \ \ \ \ \ \ \ \ 
\times\sum_{i_{41}=0}^{m_{2}}\binom{m_{2}}{i_{41}}\sum_{j_{41}=0}^{m_{1}}\sum_{j_{411}=0}^{j_{41}}\binom{m_{1}}{j_{41}}\binom{j_{41}}{j_{411}}\frac{x_1^{j_{411}}x_2^{i_{41}+j_{41}-j_{411}}x_3^{m_{2}-i_{41}+m_{3}+m_{1}-j_{41}}}{m_{1}!m_{2}!m_{3}!} .\notag 
\end{align*}
Therefore, from the definition of triple-indexed poly-Bernoulli numbers, we have
\begin{align*}
&\sum_{m_{1},m_{2},m_{3}=0}^{\infty}{\mathbb B}_{m_{1},m_{2},m_{3}}^{\star,(s_1,s_2,s_3)}\frac{x_1^{m_{1}}x_2^{m_{2}+1}x_3^{m_{3}+1}}{m_{1}!m_{2}!m_{3}!}+\sum_{m_{1},m_{2},m_{3}=0}^{\infty}{\mathbb B}_{m_{1},m_{2},m_{3}}^{\star,(s_1,s_2,s_3)}\frac{x_1^{m_{1}}x_2^{m_{2}}x_3^{m_{3}+2}}{m_{1}!m_{2}!m_{3}!} \notag \\
&=\sum_{m_{1},m_{2},m_{3}=0}^{\infty}{\mathbb B}_{m_{1},m_{2},m_{3}}^{(s_1,s_2,s_3)}\frac{x_1^{m_{1}}x_2^{m_{2}+1}x_3^{m_{3}+1}}{m_{1}!m_{2}!m_{3}!}+\sum_{m_{1},m_{2},m_{3}=0}^{\infty}{\mathbb B}_{m_{1},m_{2},m_{3}}^{(s_1,s_2,s_3)}\frac{x_1^{m_{1}}x_2^{m_{2}}x_3^{m_{3}+2}}{m_{1}!m_{2}!m_{3}!} \notag \\
& \ \ \ \ \ \ \ \ +\sum_{m_{2},m_{1},m_{3}=0}^{\infty}B_{m_{2}}{\mathbb B}_{m_{1},m_{3}}^{(s_1+s_2,s_3)}\sum_{i_{2}=0}^{m_{2}}\binom{m_{2}}{i_{2}}\sum_{j_{21}=0}^{m_{1}}\binom{m_{1}}{j_{21}}\frac{x_1^{j_{21}}x_2^{i_{2}+m_{1}-j_{21}}x_3^{m_{2}-i_{2}+m_{3}+1}}{m_{1}!m_{2}!m_{3}!} \notag \\
& \ \ \ \ \ \ \ \ +\sum_{m_{3},m_{1},m_{2}=0}^{\infty}B_{m_{3}}{\mathbb B}_{m_{1},m_{2}}^{(s_1,s_2+s_3)}\sum_{j_{32}=0}^{m_{2}}\binom{m_{2}}{j_{32}}\frac{x_1^{m_{1}}x_2^{j_{32}+1}x_3^{m_{3}+m_{2}-j_{32}}}{m_{1}!m_{2}!m_{3}!} \notag \\
& \ \ \ \ \ \ \ \ +\sum_{m_{3},m_{1},m_{2}=0}^{\infty}B_{m_{3}}{\mathbb B}_{m_{1},m_{2}}^{(s_1,s_2+s_3)}\sum_{j_{32}=0}^{m_{2}}\binom{m_{2}}{j_{32}}\frac{x_1^{m_{1}}x_2^{j_{32}}x_3^{m_{3}+m_{2}-j_{32}+1}}{m_{1}!m_{2}!m_{3}!} \notag \\
& \ \ \ \ \ \ \ \ +\sum_{m_{2},m_{3},m_{1}=0}^{\infty}B_{m_{2}}B_{m_{3}}{\mathbb B}_{m_{1}}^{(s_1+s_2+s_3)}\sum_{i_{41}=0}^{m_{2}}\binom{m_{2}}{i_{41}}\sum_{j_{41}=0}^{m_{1}}\sum_{j_{411}=0}^{j_{41}}\binom{m_{1}}{j_{41}}\binom{j_{41}}{j_{411}}\\
&\ \ \ \ \ \ \ \ \hspace{3mm}\times \frac{x_1^{j_{411}}x_2^{i_{41}+j_{41}-j_{411}}x_3^{m_{2}-i_{41}+m_{3}+m_{1}-j_{41}}}{m_{1}!m_{2}!m_{3}!} \notag \\
&=\sum_{m_{1},m_{2},m_{3}=0}^{\infty}m_2m_3{\mathbb B}_{m_{1},m_{2}-1,m_{3}-1}^{(s_1,s_2,s_3)}\frac{x_1^{m_{1}}x_2^{m_{2}}x_3^{m_{3}}}{m_{1}!m_{2}!m_{3}!}+\sum_{m_{1},m_{2},m_{3}=0}^{\infty}m_3(m_3-1){\mathbb B}_{m_{1},m_{2},m_{3}-2}^{(s_1,s_2,s_3)}\frac{x_1^{m_{1}}x_2^{m_{2}}x_3^{m_{3}}}{m_{1}!m_{2}!m_{3}!} \notag \\
& \ \ \ \ \ \ \ \ +\sum_{m_{2},m_{1},m_{3}=0}^{\infty}m_3B_{m_{2}}B_{m_{1},m_{3}-1}^{(s_1+s_2,s_3)}\sum_{i_{2}=0}^{m_{2}}\binom{m_{2}}{i_{2}}\sum_{j_{21}=0}^{m_{1}}\binom{m_{1}}{j_{21}}\frac{x_1^{j_{21}}x_2^{i_{2}+m_{1}-j_{21}}x_3^{m_{2}-i_{2}+m_{3}}}{m_{1}!m_{2}!m_{3}!} \notag \\
& \ \ \ \ \ \ \ \ +\sum_{m_{3},m_{1},m_{2}=0}^{\infty}m_2B_{m_{3}}{\mathbb B}_{m_{1},m_{2}-1}^{(s_1,s_2+s_3)}\sum_{j_{32}=0}^{m_{2}}\binom{m_{2}}{j_{32}}\frac{x_1^{m_{1}}x_2^{j_{32}}x_3^{m_{3}+m_{2}-j_{32}}}{m_{1}!m_{2}!m_{3}!} \notag \\
& \ \ \ \ \ \ \ \ +\sum_{m_{3},m_{1},m_{2}=0}^{\infty}m_3B_{m_{3}-1}{\mathbb B}_{m_{1},m_{2}}^{(s_1,s_2+s_3)}\sum_{j_{32}=0}^{m_{2}}\binom{m_{2}}{j_{32}}\frac{x_1^{m_{1}}x_2^{j_{32}}x_3^{m_{3}+m_{2}-j_{32}}}{m_{1}!m_{2}!m_{3}!} \notag \\
& \ \ \ \ \ \ \ \ +\sum_{m_{2},m_{3},m_{1}=0}^{\infty}B_{m_{2}}B_{m_{3}}{\mathbb B}_{m_{1}}^{(s_1+s_2+s_3)}\sum_{i_{41}=0}^{m_{2}}\binom{m_{2}}{i_{41}}\sum_{j_{41}=0}^{m_{1}}\sum_{j_{411}=0}^{j_{41}}\binom{m_{1}}{j_{41}}\binom{j_{41}}{j_{411}}\\
& \ \ \ \ \ \ \ \ \hspace{3mm}\times \frac{x_1^{j_{411}}x_2^{i_{41}+j_{41}-j_{411}}x_3^{m_{2}-i_{41}+m_{3}+m_{1}-j_{41}}}{m_{1}!m_{2}!m_{3}!}. \notag 
\end{align*}
Comparing the coefficients of this equation leads to Theorem \ref{Thm2}.

\begin{Remark}
Similar discussion can be made for $n$-indexed star-poly-Bernoulli numbers with $n\geq 4$.
For example, in the case of $n=4$, it is expected that single, double, triple and quadruple-indexed poly-Bernoulli numbers
will appear. Next stage is to find the regularity there and obtain the general formula.
\end{Remark}

%
%
%
%
%
%

%
%

%
\bigskip
\noindent
\textsc{Yuna Baba}\\
Faculty of Science and Technology, \\
 Sophia University, \\
 7-1 Kio-cho, Chiyoda-ku, Tokyo, 102-8554, Japan\\

\medskip

\noindent
\textsc{Maki Nakasuji}\\
Department of Information and Communication Science, Faculty of Science  and Technology, \\
 Sophia University, \\
 7-1 Kio-cho, Chiyoda-ku, Tokyo, 102-8554, Japan \\
 \texttt{nakasuji@sophia.ac.jp}\\
and\\
Mathematical Institute, \\
Tohoku University, \\
Sendai 980-8578, Japan \\
\medskip

\noindent
\textsc{Mika Sakata}\\
Department of Health and Sport management, School of Physical Education,\\
Osaka University of Health and Sport Sciences, \\
 1-1 Asashirodai, Kumatori-cho, Sennan-gun, 
Osaka, 590-0496, Japan \\
 \texttt{m.sakata@ouhs.ac.jp}\\

\end{document}